\documentclass{article}
\usepackage[utf8]{inputenc}
\def\MR#1{\href{http://www.ams.org/mathscinet-getitem?mr=#1}{MR#1}}
\usepackage{amssymb}
\usepackage{hyperref}
\usepackage{amsmath}
\newtheorem{lem}{Lemma}
\newtheorem{ter}{Theorem}

\newtheorem{utv}{Proposition}
\newtheorem{sled}{Corollary}
\title{
The $\alpha$-representation for the characteristic function of a matroid}
\author{Eduard~Yu.~Lerner}
\date{}
\begin{document}
\maketitle
\begin{abstract}
Let $M=(E,\mathcal B)$ be an $\mathbb F_q$-linear matroid; denote by ${\mathcal B}$ the family of its bases, $s(M;\alpha)=\sum_{B\in\mathcal B}\prod_{e \in B} \alpha_e$, where ${\alpha_e\in \mathbb F_q}$. According to the Kontsevich conjecture stated in 1997, the number of nonzero values of $s(M;\alpha)$ is a polynomial with respect to $q$ for all matroids. This conjecture was disproved by P.~Brosnan and P.~Belkale.
In this paper we express the characteristic polynomial of the dual matroid $M^\perp$
in terms of the ``correct'' Kontsevich formula (for $\mathbb F_q$-linear matroids). This repre\-sen\-ta\-tion generalizes the formula for a flow polynomial of a graph which was obtained by us earlier (and with the help of another technique). In addition, generalizing the correlation (announced by us earlier) that connects flow and chromatic polynomials, we define the characteristic polynomial of $M^\perp$ in two ways, namely, in terms of characteristic polynomials of~$M/A$ and $M|_A$, respectively, $A\subseteq E$. The latter expressi\-ons are close to convolution-multiplication formulas established by V.~Reiner and J.~P.~S.~Kung.
\end{abstract}


\section{Introduction. Statements of main results}

Let us use the following denotations: $\mathbb F_q$ is a finite field of an odd characteristic~$p$; $M$ is an $\mathbb F_q$-linear matroid on a set~$E$;
${\,\tt M\,}$ is the matrix with elements in the field~$\mathbb F_q$ that defines the mentioned matroid; we enumerate columns of this matrix with elements of the set~$E$ and do rows with elements of the set~$V$. We understand a matroid $M$ as a pair $(E,r)$, where $r$ is the rank function defined on subsets~$A$ of the set~$E$. In the case of an $\mathbb F_q$-linear matroid, $r(A)$ equals the rank of the set of columns of the matrix~${\,\tt M\,}$ with indices in~$A$
(see, for example, \cite{Oxley,aigner,reiner1} for more details).

In what follows for convenience we assume that $r(E)=|V|$ (all linearly dependent rows are already deleted from the matrix ${\,\tt M\,}$). For example, if $G$ is a connected multigraph and $M$ is its cyclic matroid, then instead of the matrix~${\,\tt M\,}$ we can consider its incidence matrix with one deleted row (chosen arbitrarily). Recall that the incidence matrix $(\varepsilon_{v e})_{v\in V,e\in E}$ obeys the formula
$$
\varepsilon_{v e}  = \left\{
  \begin{array}{ll}
   -1, & \mbox{if $i(e)=v$},  \\
   1, & \mbox{if $f(e)=v$},  \\
   0, & \mbox{if $e$ is nonincident to $v$},
  \end{array} \right.,
$$
where $i(e)$ ($f(e)$) is the origin (the endpoint) of the edge~$e$ different from a loop;
in the case of a loop~$e$ the corresponding column of the matrix~$\varepsilon$ consists only of zeros.
Denote by $M^\perp$ the dual to~$M$ matroid $(E,r^\perp)$: $r^\perp(A)=|A|-r(A)+r(E\setminus A)$ for any $A: A\subseteq E$. Recall that the characteristic polynomial $\chi_M$ of the matroid~$M$ obeys the formula
\begin{equation}
\label{chiM}
\chi_M(x)=\sum_{A:A\subseteq E} (-1)^{|A|} x^{r(E)-r(A)}.
\end{equation}  

Let us associate each element $e$ of the set $E$ with a variable $\alpha_e$ that takes on values in~$\mathbb F^*_q$, where $\mathbb F^*_q$ is the totality of nonzero elements of the field~$\mathbb F_q$. The $\alpha$-representation defines the expression for $\chi_{M^\perp}(q)$ in terms of a linear combination of Legendre symbols of some polynomial $s({\,\tt M\,};\alpha)$ of the variables~$\alpha$.
Let us define this polynomial.

Let $\mathcal B({\,\tt M\,})$ be the family of all possible bases of the matroid generated by the matrix~${\,\tt M\,}$. By definition,
$\mathcal B({\,\tt M\,})$ is the family of maximum cardinality sets~$B$ of indices from the set~$E$ such that the corresponding columns of the matrix ${\,\tt M\,}$ are linearly independent. Denote by ${\,\tt M\,}|_B$ the nondegenerate square submatrix of the matrix~${\,\tt M\,}$ formed by columns whose numbers belong to the set~$B$ (and all rows of the matrix~${\,\tt M\,}$). Put
$$s({\,\tt M\,};\alpha)=\sum_{B\in{\mathcal B}({\,\tt M\,})}{\det}^2 ({\,\tt M\,}|_B) \prod_{e \in B} \alpha_e.$$
For the matroid~$M$ of zero rank (it is formed by the matrix whose row set is empty and column set $E$ consists of loops of the matroid~$M$) we put $s({\,\tt M\,};\alpha)\equiv 1$.

Note that in the case of a cyclic matroid of a connected multigraph any basis $B$ is a spanning tree, therefore for the matrix ${\,\tt M\,}$ (which is obtained from the matrix~$\varepsilon$ by deleting an arbitrary row) the submatrix ${\,\tt M\,}|_B$ is (up to the sign of elements) a permutation matrix, consequently, $\det^2 ({\,\tt M\,}|_B)=1$. In this case the sum $s({\tt M};\alpha)$ depends only on bases of the matroid~$M$ and equals
\begin{equation}
\sum_{B\in\mathcal B}\prod_{e \in B} \alpha_e.
\label{simplesum}
\end{equation}
In a general case an $\mathbb F_q$-linear matroid does not necessarily has a representation in terms of the matrix~${\,\tt M\,}$ such that $s({\,\tt M\,};\alpha)$ coincides with sum~(\ref{simplesum})(in what follows we illustrate this property by considering the matroid~$U_{2,4}$ as an example).

Let~$W\subseteq V$. Denote by ${\,\tt M\,}/W$ the matrix obtained from the matrix ${\,\tt M\,}$ by deleting rows whose numbers belong to the set~$W$; this matrix corresponds to some matroid. Evidently, the rank of this matroid equals $|V|-|W|$. 

Denote by $W^*_{{\,\tt M\,}}(\alpha)$ an arbitrary minimum cardinality subset $W$ of the set~$V$, for which the sum $s({\,\tt M\,}/W;\alpha)$ differs from zero. Denote by $r^*({\,\tt M\,};\alpha)$ the difference $|V|-|W^*_{{\,\tt M\,}}(\alpha)|$.

Recall that we consider the field~$\mathbb F_q$ of an odd characteristic, i.e., $q=p^d$, where $p$ is an odd prime, and $d\in\mathbb N$. Denote by $\eta$ the multiplicative quadratic character of the field~$\mathbb{F}_q$: $\eta(0)=0$, in other cases $\eta(x)=1$ or $\eta(x)=-1$ depending on whether $x$ is a square in the field $\mathbb F_q$ or not. For $d=1$ the function $\eta$ coincides with the Legendre symbol of the residue field modulo prime~$p$. Finally, let us define a function $g(q,n)$, where $q$ is the number of elements in the mentioned field, $n\in\mathbb N$, by the formula
\[
g(q,n) = \left\{ \begin{array}{ll}
   1/q^{n/2},& \textrm{if}~p \bmod 4= 1,\,  n\bmod 2= 0,  \\
   1/(-q)^{n/2},&  \textrm{if}~p \bmod 4=3,\,  n\bmod 2= 0,\\
   0,& \textrm{if}~ n\bmod 2=1.
     \end{array} \right.
\]

\begin{ter}{\rm \textbf{(The main theorem)}}\
\label{ter:1}
Let us use the following denotations: $\mathbb F_q$ is a finite field of an odd characteristic, $M$ is an $\mathbb F_q$-linear matroid, and ${\,\tt M\,}$ is some matrix of its representation. The following formula is valid:
\begin{equation} \label{eq:3}
\chi_{M^\perp} (q) =
\sum_{\alpha \in (\mathbb F^*_q)^E}
g(q,r^*({\,\tt M\,};\alpha))\,
\eta(s({\,\tt M\,}/W^*_{{\,\tt M\,}}(\alpha);\alpha)).
\end{equation}
\end{ter}

In December 1997, when giving a talk at the Gelfand seminar in Rutgers University, Maxim Kontsevich proposed a conjecture that for any matroid~$M$ the number of nonzero values of~(\ref{simplesum}) for $\alpha  \in \mathbb{F}_q^E$ is a polynomial with respect to~$q$. Though the conjecture was never published, it has aroused the interest of experts in combinatorics (see~\cite{stanleyArticle,chung}). Against expectations, sometime later this conjecture was refuted~\cite{belk}. Formula~(\ref{eq:3}) represents a ``correct'' variant of the Kontsevich conjecture for arbitrary $\mathbb{F}_q$-linear matroids.

In the case of a graphic matroid, $\chi_{M^\perp}$ is a flow polynomial of the graph, and $s({\,\tt M\,};\alpha)$ is
an arbitrary cofactor of the weighted Laplacian matrix of the graph.
This representation was obtained by us earlier~\cite{ourArX}. Note that methods that use the Fourier transformation described in~\cite{ourArX} may appear to be inapplicable in a general case. Moreover, formula~(\ref{simplesum}), which is used in~\cite{ourArX} for~$s(M;\alpha)$, in a general case is incorrect. Namely, for a matroid~$M$ the existence of a matrix~${\,\tt M\,}$ such that $s({\,\tt M\,};\alpha)$ coincides with~(\ref{simplesum}) is not guaranteed.

Consider, as an example, the matroid $U_{2,4}$ and the corresponding to it matrix~${\,\tt M\,}$. Since we need to ensure the condition ${\det}^2 ({\,\tt M\,}|_B)=1$ for the basis~$B$ consisting of columns $\{1,2\}$ of this matrix, the linear transformation with the determinant $\pm 1$ turns the first two columns to the form ${1}\choose{\,0\,}$ and~${\,0\,}\choose{1}$, respectively.
Below we assume that columns $\{1,2\}$ take these forms. Then the same requirement that ${\det}^2 ({\,\tt M\,}|_B)=1$ but for each basis~$B$ consisting of one of columns $\{1,2\}$ and one of those $\{3,4\}$ is equivalent to the requirement that the last two columns should consist of elements~{$\pm 1$}. The nonzero determinant of the matrix consisting of columns $\{3,4\}$ can be equal only to~$\pm 2$. We get
$$
s({\,\tt M\,};\alpha)=\alpha_1\alpha_2+\alpha_1\alpha_3+\alpha_1\alpha_4+\alpha_2\alpha_3+\alpha_2\alpha_4+4\alpha_3\alpha_4,
$$
which differs from formula~(\ref{simplesum}).

Let us illustrate the main correlation~(\ref{eq:3}) in the case of the matroid~$U_{2,4}$. Note that the matroid~$U_{2,4}$ is dual to itself, while the left-hand side of equality~(\ref{eq:3}) is
$\chi_{\,U_{2,4}^\perp}(q)=\chi_{\,U_{2,4}}(q)=(q-1)(q-3)$. If~$U_{2,4}$ is defined by the matrix
\begin{equation}
\label{MM}
{\,\tt M\,}=\left(\begin{matrix}1&0&1&1\\0&1&1&-1\end{matrix}\right),
\end{equation}
then we can prove that $r^*({\,\tt M\,};\alpha)=0$ only in $q-1$ cases (namely, when ${\alpha_1=\alpha_2}$, $\alpha_3=\alpha_4$, and $\alpha_1=-2\alpha_3$). Therefore, for~$U_{2,4}$ with matrix repre\-sen\-ta\-tion~(\ref{MM}) Theorem~\ref{ter:1} is equivalent to the equality
$$
(q-1)(q-4)=g(q,2)\sum_{\alpha \in (\mathbb F^*_q)^4}
\eta(\alpha_1\alpha_2+\alpha_1\alpha_3+\alpha_1\alpha_4+\alpha_2\alpha_3+\alpha_2\alpha_4+4\alpha_3\alpha_4).
$$

Another theorem which is also proved in this paper calculates $\chi_{M^\perp} (q)$ for an arbitrary matroid~$M$
(which is not necessarily representable over some field). The right-hand sides of identities given below equal sums of polynomials, i.e., normalized characteristic functions of matroids defined on subsets of the set~$E$. Therefore these identities (as distinct from
Theorem~\ref{ter:1}) impose no constraints on values of the argument~$q$.

Let $A\subseteq E$. Denote by $M|_A$ and $M/A$ matroids obtained by restricting $M$ to set~$A$ 
and contracting this set, respectively
(see details, for example, in~\cite{Oxley,aigner,reiner1}).

\begin{ter}
Let $M=(E,r)$ be an arbitrary matroid on the set~$E$ with the rank function~$r$. The following formulas are valid:
\label{ter:2}
\begin{eqnarray}
\label{one}
& \chi_{M^\perp}(q)
=  (q-1)^{|E|} \sum\limits_{A:A\subseteq E} \left(\frac{q}{1-q}\right)^{|A|}
\frac{\chi_{M|_A}(q)}{q^{r(A)}},&\\
\label{two}
&\chi_{M^\perp}(q)
=  q^{-r(E)} \sum\limits_{A:A\subseteq E} (-1)^{|E|-|A|}(q-1)^{|A|}\chi_{M/A}(q).&
\end{eqnarray}
\end{ter}
In the statement of identity~(\ref{one}) one can use the function ${\zeta_q(z)=\frac1{1-q^{-z}}}$; more precisely, values of this function at points~1 and -1, i.e., $\zeta_q(1)=\frac1{1-q^{-1}}$, $\zeta_q(-1)=\frac1{1-q}$. Then we can rewrite identity~(\ref{one}) as follows:
\begin{equation}
\label{onezeta}
\chi_{M^\perp}(q)\,(\zeta_q(-1))^{|E|}
=\sum_{A:A\subseteq E}
(-1)^{|E|-|A|}\,
\frac{\chi_{M|_A}(q)}{q^{r(A)}}\,
(\zeta_q(1))^{|A|}.
\end{equation}

With the help of the inclusion-exclusion formula one can easily prove the equivalence of correlations~(\ref{one}) and~(\ref{two}). Really, as is well known, $\left( M|_A \right)^\perp=M^\perp/(E\setminus A)$. Applying the Mobius inversion formula (in the case under consideration it coincides with the inclusion-exclusion formula) to identity~(\ref{onezeta}), we can rewrite the latter as follows:
\begin{equation}
\label{twozeta}
\frac{\chi_M(q)}{q^{r(E)}}(\zeta_q(1))^{|E|}=\sum_{A:A\subseteq E} (\zeta_q(-1))^{|A|}\,\chi_{M^\perp/(E\setminus A)}(q).
\end{equation}
Replacing the matroid~$M$ in formula~(\ref{twozeta}) with $M^\perp$ (recall that $\left(M^\perp\right)^\perp=M$) and performing trivial transformations, we get formula~(\ref{two}). Therefore, correlations~(\ref{one}) and~(\ref{two}) are equivalent, so it suffices to prove only one of them.

In the case when $M$ is a cyclic matroid of some graph~$G$, $\chi_{M^\perp}(q)$ coincides with the flow polynomial of the graph $F_G(q)$, and $\chi_M(q)$ is expressed in terms of the chromatic polynomial~$P_G(q)$ (see~\cite{tutte-book,aigner,reiner1}). In this case identities~(\ref{one})--(\ref{twozeta}) express the interconnection between flow and chromatic polynomials.

In particular, in this case identity~(\ref{onezeta}) means that
\begin{equation}
\label{kuptscov1}
F_G (q)\left( {\zeta _q ( - 1)} \right)^{|E(G)|}  =
\sum\limits_{H:H \subseteq G}  ( - 1)^{|E(G)| -|E(H)|}
\frac{P_H (q)}{q^{|V(H)|}}\left( {\zeta _q (1)} \right)^{|E(H)|},
\end{equation}
where the sum is taken over all subgraphs~$H$ of the graph~$G$, while $V(G)$ is the set of vertices of the graph~$G$ (note that in formula~(\ref{kuptscov1}) $V(H)=V(G)$) and $E(G)$ is the set of its edges. In this case identity~(\ref{twozeta}) is equivalent to an even more simple correlation, namely,
\begin{equation}
\label{kuptscov2}
\frac{P_G (q)}{q^{|V(G)|}}\left(\zeta _q (1)\right)^{|E(G)|}
=\sum_{H: H \subseteq G} \left(\zeta _q(-1)\right)^{|E(H)|} F_H (q).
\end{equation}

Finally, in the case of a cyclic matroid of the graph~$G$ identity~(\ref{two}) takes the form
\begin{equation}
\label{twograph}
F_G(q) q^{|V(G)|}=\sum_{H: H\subseteq G} (-1)^{|E(G)|-E(H)|}(q-1)^{|E(H)|}P_{G/H}(q),
\end{equation}
where $G/H$ is the contracted graph~$G$. Note that $P_{G/H}(q)\equiv 0$, provided that the graph $G/H$ has loops. For example, in the case of the graph~$C_4$, the right-hand side of~(\ref{twograph}) contains the following terms: the value $P_{C_4}(q)$ that corresponds to the case of $E(H)=\emptyset$; 4 terms $-(q-1)P_{K_3}(q)$ that correspond to the case of $|E(H)|=1$; 6 terms $(q-1)^2 P_{K_2}(q)$ that correspond to the case of $|E(H)|=2$; and one term $(q-1)^4 q$, ($H=C_4$); there are no terms that correspond to the contraction of the graph consisting of three edges, because in this case there appears a loop in the graph $G/H$. Therefore, for the graph~$C_4$ identity~(\ref{twograph}) takes the form
$$
(q-1)q^4=P_{C_4}(q)-4(q-1)P_{K_3}(q)+6(q-1)^2 P_{K_2}(q)+(q-1)^4 q;
$$
one can be verify this fact immediately.

Equalities~(\ref{one},\ref{two}) are close to convolution formulas given in papers by W.~Kook, V.~Reiner, and D.~Stanton and those by V.~Reiner and J.~P.~S.~Kung (\cite{reiner3,reiner2,kung1,kung2}). Moreover, one can etablish formulas~(\ref{one},\ref{two}) (with the help of certain transformations) as particular cases of identities proved in~\cite{kung2}. One can also prove them immediately by using Definition~(\ref{chiM}). We give them in this paper, because the existence of these formulas in the case of the flow polynomial (as well as the existence of the $\alpha$-representation for it) evidently follows from the general theory of Feynman amplitudes. In Section~3 we establish the connection between the flow and chromatic polynomials within the framework of the theory of Feynman amplitudes. In Section~4 we prove two other variants of Theorem~\ref{ter:2}, which differ in the generality of statements, and discuss the connection of this theorem with those obtained in earlier papers. Section~2 is dedicated to the proof of Theorem~\ref{ter:1}. In conclusion we summarize the obtained results and discuss further perspectives of the research in the application of methods developed within the framework of the theory of Feynman amplitudes for flow polynomials.

For structuring this paper we distinguish between theorems, lemmas, propositions, and corollaries in it. We formulate the most important results as theorems. The destination of lemmas and corollaries is evident. Propositions play a special role. Namely, we state well-known facts as propositions; we prove them when failing in finding proper references.

\section {Proof of Theorem~\ref{ter:1}}

The proof of Theorem~\ref{ter:1} essentially uses an analog of the Laplacian matrix of the graph (more precisely, an analog of the submatrix obtained from this matrix by deleting the last row and the last column). This is the following square matrix~$L({\,\tt M\,};\alpha)$ of the order~$r(M)$:
$$
L({\,\tt M\,};\alpha)={\,\tt M\,} \Lambda {\,\tt M\,}^T,
$$
where $\Lambda$ is the $|E|\times|E|$-diagonal matrix whose diagonal elements equal $\alpha_e$,
$e\in E$, and $T$ is the transposition sign.

\begin{lem}
\label{lem:1}
Let ${\mathbf x}=(x_v, v\in V)$ be a row of variables that take on values in the field~$\mathbb F_q$; let $Q({\mathbf x},\alpha)={\mathbf x} L({\,\tt M\,};\alpha){\mathbf x}^T$ be the quadratic form with the matrix $L({\,\tt M\,};\alpha)$; for any $b\in\mathbb F_q$ we denote by $N_b$ the number of solutions of the equation
\begin{equation}
\label{eq}
Q({\mathbf x},\alpha)=b
\end{equation}
with respect to variables ${\mathbf x}\in \mathbb F_q^V$ and $\alpha\in \left(F^*_q\right)^E$. Then\\
1) $N_b=N_1$ for any $b\in F^*_q$;\\
2) $\chi_{M^\perp}(q)=(N_0-N_1)/q^{|V|}$.
\end{lem}
\textbf{Proof:} The first proposition is almost evident. Really, let $b$ be a fixed element of $\mathbb F_q$. Then any collection $({\mathbf x},\alpha)$, ${\mathbf x}\in\mathbb F_q^V$, $\alpha\in\left(\mathbb F^*_q\right)^E$, such that $Q({\mathbf x},\alpha)=1$ bijectively corresponds to the collection $({\mathbf x},\beta)$, where $\beta=b \alpha$. Evidently, $Q({\mathbf x},\beta)=b$, i.e., we have established a bijection between solutions of Eqs.~(\ref{eq}) with various $b\in\mathbb F_q^*$.

The second proposition is less trivial; by certain modifications one can reduce it to fundamental results obtained by H.~Crapo and G.-C.~Rota. Assume that ${\mathbf y}({\mathbf x})={\mathbf x} {\,\tt M\,}$, i.e., ${\mathbf y}$ is the row $(y_e,e\in E)$. Evidently, $Q({\mathbf x},\alpha)=\sum_{e:e\in E}\alpha_e y_e^2$.

Denote by $Q'({\mathbf x},\alpha)$ the value of the function $\sum_{e:e\in E}\alpha_e y_e$; let $N'_b$, $b\in\mathbb F_q$, be the number of solutions $({\mathbf x},\alpha)$ of the equation $Q'({\mathbf x},\alpha)=b$. Let us prove that $N'_b=N_b$ for any $b\in\mathbb F_q$.

Really, let us associate the collection $({\mathbf x},\alpha)$ with that $({\mathbf x},\gamma)$ such that $\gamma_e=\alpha_e y_e$, if $y_e\neq 0$, otherwise $\gamma_e=\alpha_e$. Then $Q({\mathbf x},\alpha)\equiv Q'({\mathbf x},\gamma)$ and therefore $N'_b=N_b$.

Let us now prove that $(N'_1-N'_0)/q^{|V|}$ coincides with the number of vectors~$\alpha$ such that ${\mathbf c\,}(\alpha)=0$, where ${\mathbf c\,}(\alpha)={\,\tt M\,}\alpha$ is the vector column $(c_v,v\in V)$.

The equation $Q'({\mathbf x},\alpha)=b$ takes the form $\sum_{v:v\in V} c_v x_v=b$.

If ${\mathbf c\,}(\alpha)=0$, then $Q'({\mathbf x},\alpha)=0$ independently of the choice of the vector ${\mathbf x}$. Note that the total number of vectors~${\mathbf x}$ equals $q^{|V|}$.

However if ${\mathbf c\,}(\alpha)\ne 0$, then the equation $Q'({\mathbf x},\alpha)=b$ is equivalent to that $c_{v'}x_{v'}=b-\sum_{v:v\ne v'} c_v x_v$, where $c_{v'}$ is the first nonzero element of the vector~${\mathbf c\,}$. The number of solutions of this equation with respect to ${\mathbf x}$ equals $q^{|V|-1}$ independently of the choice of $b\in\mathbb F_q$, because by fixing arbitrary $(x_v, v\in V \setminus \{ v'\})$ we uniquely define~$x_{v'}$.

Thus, we have proved that $(N_1-N_0)/q^{|V|}$ coincides with the number of vectors~$\alpha$ such that ${\mathbf c\,}(\alpha)=0$, i.e., with the cardinality of the everywhere nonzero kernel of the map ${\,\tt M\,}\alpha$ (recall that $\alpha_e\in\mathbb F_q^*$). If ${\,\tt M\,}$ is the incidence matrix of the graph~$G$, then the equation ${\,\tt M\,}\alpha=0$ guarantees (by definition) that the collection $(\alpha_e,e\in E(G))$ is a flow on edges of the directed graph~$G$. The number of everywhere nonzero flows equals the value of the flow poly\-no\-mial~$F_G(q)$. In a general case, as was proved by H.~Crapo and G.-C.~Rota (in this connection see, for example, \cite{reiner2}), cardinality of the everywhere nonzero kernel of the map ${\,\tt M\,}\alpha$ coincides with the value~$\chi_{M^\perp}(q)$. $\square$

Note that the following proposition is an evident generalization of Lemma~\ref{lem:1}.

\begin{lem}
\label{utv:1}
Let $j\in \mathbb N$. Put $Q_m({\mathbf x},\alpha)=\sum_{e:e\in E} \alpha_e (y_e)^j$, where ${\mathbf y}\equiv {\mathbf y}({\mathbf x})={\mathbf x} {\,\tt M\,}$. Then for any $b\in\mathbb  F_q$ the number $N_b(j)$ of solutions of the equation ${Q_j({\mathbf x},\alpha)=b}$ is independent of~$j$. Moreover,
$$\chi_{M^\perp}(q)=(N_0(j)-N_1(j))/q^{|V|}.$$
\end{lem}

Really, the proof of Lemma~\ref{lem:1} is based on the proof of the coincidence of $N_b(2)$ (in denotations of the lemma, $N_b$) and $N_b(1)$ (in denotations of the lemma, $N'_b$). This coincidence was proved with the help of the equality $Q_2({\mathbf x},\alpha)=Q_1({\mathbf x},\gamma)$, where there is a bijective correspondence between $\gamma$ and $\alpha$.
One can easily see that $Q_{j+1}({\mathbf x},\alpha)=Q_j({\mathbf x},\gamma)$ with any $j\in\mathbb N$, which proves Lemma~\ref{utv:1}.

\begin{utv}
\label{lem:2}
Let ${\mathbf x}=(x_1,\ldots,x_n)$ be a row of variables whose values belong to the field~$\mathbb F_q$; let $Q_B({\mathbf x})={\mathbf x} B {\mathbf x}^T$ be a quadratic form with a fixed symmetric $n\times n$-matrix~$B$ with elements in the mentioned field. Let $m$ be the rank of the matrix~$B$ and let $d$ be its arbitrary nonzero principal minor of the order~$m$. The number~$N$ of solutions of the equation $Q_B({\mathbf x})=0$ with respect to variables ${\mathbf x}\in \mathbb F_q^n$ obeys the formula
$$
N=\left\{\begin{array}{ll}
q^{n-m}\left(q^{m-1}+(q-1)q^{\frac{m-2}{2}}\eta((-1)^{m/2} d) \right), & \textrm{if} ~ m\bmod 2= 0,\\
      q^{n-1},& \textrm{if} ~ m\bmod 2=1.
     \end{array} \right.
$$
\end{utv}

Proposition~\ref{lem:2}, apparently, has been known since works of C.~Chevalley. The case when~$m=n$ coincides with theorem~2.E in (\cite{schmidt}, Ch.~4); see also \cite[Theorem 6.26]{lidl}). The general case can be easily reduced to the case of $m=n$.

\textbf{Proof of Theorem~\ref{ter:1}}: Let us first state Lemma~\ref{lem:1} in terms of zeros of quadratic forms mentioned in Proposition~\ref{lem:2}. Graduating $N_0$ with respect to all possible values of~$\alpha$, {we rewrite the main equality of this lemma as follows:}
\begin{equation}
\label{Cheval}
q^{|V|}\chi_{M^\perp}(q)=
\sum_{\alpha\in \left(\mathbb F_q^*\right)^E}
\left(N(\alpha)-\frac{|\mathbb F_q^V|-N(\alpha)}{q-1}\right),
\end{equation}
where $N(\alpha)$ is the number of solutions of the equation $Q({\mathbf x},\alpha)=0$ with respect to~${\mathbf x}$ with fixed~$\alpha$, while~$|\mathbb F_q^V|=q^{|V|}$.

For calculating~$N(\alpha)$ let us use Proposition~\ref{lem:2}. Note that if $\mathop{\rm rank} L({\,\tt M\,};\alpha)$ is odd, then the corresponding term in sum~(\ref{Cheval}) equals zero, because in this case it holds that
$$
N(\alpha)=\frac{q^{|V|}-N(\alpha)}{q-1}=q^{|V|-1}.
$$
If $m=\mathop{\rm rank} L({\,\tt M\,};\alpha)$ is even, then we get
$$
N(\alpha)-\frac{q^{|V|}-N(\alpha)}{q-1}=q^{|V|-m}\left(q q^{\frac{m-2}{2}}\eta((-1)^{m/2} d(\alpha)) \right),
$$
where $d(\alpha)$ is an arbitrary nonzero principal minor of the matrix $L({\,\tt M\,};\alpha)$ of the order~$m$. By elementary transformations we rewrite formula~(\ref{Cheval}) as follows:
$$
\chi_{M^\perp} (q) =
\sum_{\alpha \in (\mathbb F^*_q)^E}
g(q,\mathop{\rm rank} L({\,\tt M\,};\alpha))\,\eta(d(\alpha)).
$$
It remains to prove that $\mathop{\rm rank} L({\,\tt M\,};\alpha)=r^*({\,\tt M\,};\alpha)$ and $d(\alpha)=s({\,\tt M\,}/W^*_{{\,\tt M\,}}(\alpha);\alpha)$.

To this end it suffices to prove that if the principal minor of the matrix $L({\,\tt M\,};\alpha)$ is obtained by deleting rows and columns whose numbers belong to the set~$W$ then this minor equals $s({\,\tt M\,}/W;\alpha)$. Note that by definition the desired minor coincides with $\det(({\,\tt M\,}/W) \Lambda ({\,\tt M\,}/W)^T)$, therefore it remains to prove the equality $\det(L({\,\tt M\,}';\alpha))=s({\,\tt M\,}',\alpha)$. The latter formula is a generalization of the matrix tree theorem; it is rather well-known (see \cite[Exercise~11 (c)]{reiner1}).
$\square$

\section{Chromatic and flow polynomials as Feynman amplitudes and their interconnection}
In this section we prove that formulas~(\ref{kuptscov1},\ref{kuptscov2},\ref{twograph}) follow from the general theory of Feynman amplitudes (FA).

Let $G=(V,E)$ be an arbitrary multigraph with connected compo\-n\-e\-nts~$K(G)$. Associate each edge~$e$ ($e\in E$) with a complex-valued function (a propagator) $\Delta_e$ of the argument $x\in\mathbb F_q$. In the usual theory of FA, $x\in\mathbb R^4$, below we consider a generalization of all definitions to the case of finite fields. In the real-valued case one considers even functions $\Delta_e$, so do we.

Associate each vertex $v$ of the graph~$G$ ($v\in V$) with a variable $x_v$, $x_v\in\mathbb F_q$. We understand a vacuum FA in a coordinate space as the sum
\begin{equation}
\label{defFG}
{\mathcal F}_G =\sum_{{\mathbf x} \in F_q^V} \prod_{e\in E} \Delta_e(x_{i(e)}-x_{f(e)}),
\end{equation}
where $i(e)$ ($f(e)$) is the origin (endpoint) of the edge~$e$.

Note that since the function~$\Delta$ is even, the value ${\mathcal F}_G$ is independent of the orientation of graph edges. Definition~(\ref{defFG}) is a calque of the case $x\in \mathbb R^4$, except the fact that in place of the sum in the real-valued case there is a multidimensional integral, and the integration is performed in all variables~$x_v$, {except} those ones which are chosen arbitrarily so that each connected component $K(G)$ contains exactly one such variable. An analogous exclusion of variables from the sum in~(\ref{defFG}), evidently, leads to the sum ${{\mathcal F}_G}'$ which differs from ${\mathcal F}_G$ by the trivial multiplier ${{\mathcal F}_G}'=q^{-|K(G)|}{\mathcal F}_G $.

Let us introduce on the field $\mathbb F_q$ a trivial norm $||x||$ which equals $1$ with all $x\in\mathbb F^*_q$. Note that in the case of the finite field $||x||=1-\delta(x)$, where $\delta(\cdot)$ is a finite delta function which is an indicator of the value $x=0$. One can easily see that if $\Delta_e\equiv||\cdot||$, then by definition ${\mathcal F}_G$ coincides with the value of the chromatic polynomial~$P_G(q)$.

Let us now associate each edge~$e$ of the graph with a variable $k_e$, ${k_e\in\mathbb F_q}$, and, generally speaking, some other propagator
 $\widetilde\Delta$. Put ${\mathbf k}=(k_e, e\in E)$. Recall that the symbol $(\varepsilon_{v e})_{v\in V,e\in E}$ denotes the incidence matrix. We under\-stand a vacuum FA in a momentum space as the sum
\begin{equation}
\label{defwFG}
\widetilde{\mathcal F}_G = \sum_{{\mathbf k}\in\mathbb F_q^E}\
 \prod_{e \in E} \widetilde\Delta_e(k_e) \prod_{v \in V} \delta \left( \sum_{e\in E(G)}
  \varepsilon_{v e}  k_e \right).
\end{equation}
This is the exact calque of the definition given for the real case except the fact that the integration in $\left(\mathbb R^4\right)^E$ is replaced with the summation over~$\mathbb F_q^E$. One can easily see that if $\widetilde\Delta_e\equiv||\cdot||$ then by definition $\widetilde{\mathcal F}_G$ coincides with the value of the flow polynomial~$F_G(q)$.

For any function $f(x)$ whose argument $x$ takes on values in~$\mathbb F_q$ we define the \emph{Fourier transformation} $\widehat{f}(k)$, $k\in \mathbb F_q$, as
$$
\widehat{f}(k)=\sum_{x\in
\mathbb F_q}{f(x)\chi( k x)},
$$
where the canonical additive character of the $\chi$
\footnote{Note
that the traditional denotation~$\chi$ used for the additive character has nothing in common with the denotation~$\chi_M$ used for the characteristic function of a matroid.}
field~$\mathbb F_q$ takes the form
\begin{equation}
\label{simplechi}
\chi(x)=\exp{(2\pi i\, \operatorname{Tr}(x)/p)},\quad \mbox{while $\operatorname{Tr}(x)=x+x^p+x^{p^2}\ldots+x^{p^{d-1}}$}.
\end{equation}
The function $f$ can be restored from the function~$\widehat f$ with the help of the formula
$$
f(x)=\frac1q\sum_{k\in \mathbb F_q}\widehat f(k)\chi(- k x)
$$
(the inverse Fourier transformation formula). Note that formulas used in the real-valued case {$x\in\mathbb R^4$} are nearly the same, except the fact that the summation operation is replaced with the integration one, and normalizing constants take on other values.

It is well known that in the real-valued case the equality $\Delta_e=\widehat{\,\widetilde{\Delta_e\!}\,}$ implies the correlation
$\widetilde{\mathcal F}_G={{\mathcal F}_G}'/c^{|V|-|K(G)|}$, where $c$ is the normalizing constant of the Fourier transformation. In the case of a finite field, under the same assumptions about $\widetilde\Delta_e$ and $\Delta_e$ one can analogously prove the same correlation, setting the normalizing constant~$c$ to~$q$, i.e.,
$q^{|V|}\widetilde{\mathcal F}_G={\mathcal F}_G$.
We are going to study a particular case of this {correlation}, when
\begin{equation}
\label{Delta}
\Delta_e(x)\equiv a+b\,\delta(x),\quad  \widetilde\Delta_e(k)\equiv \tilde a+\tilde b\,\delta(k).
\end{equation}
and represent it not only as a corollary of the theory of FA discussed above but also as a well-known combinatorial identity.
Note that the Fourier trans\-for\-ma\-ti\-on reduces propagators~(\ref{Delta}) to the form $\widehat{a+b\,\delta}=b+aq\,\delta$; {this defines the ``correctness'' (from the point of view of the theory of FA) of the statement of the following proposition.}

\begin{utv}
\label{ter:3}
Assume that propagators of FA in coordinate and momentary spaces obey formulas~(\ref{Delta}), and vacuum amplitudes ${\mathcal F}_G(q;a,b)$ and $\widetilde{\mathcal F}_G(q;\tilde a,\tilde b)$ do formulas~(\ref{defFG}) and~(\ref{defwFG}), correspondingly. Then
\begin{equation}
\label{eq:terfourier}
q^{|V|}\widetilde{\mathcal F}_G(q;a,b)={\mathcal F}_G(q;b,aq).
\end{equation}
\end{utv}

We state equality~(\ref{eq:terfourier}) as a proposition rather than a theorem, because, as follows from considerations given below, the equality equivalent to it has been known in combinatorics as a correlation between the so-called ``bad flow polynomial'' and ``bad coloring polynomial''.

Before we proceed to the mentioned considerations, let us write down important recurrent correlations for ${\mathcal F}_G(q;a,b)$ and $\widetilde{\mathcal F}_G(q;\tilde a,\tilde b)$.

Fix some edge~$e$ of the graph~$G$. Let us denote by $G'_e$ the graph obtained from~$G$ by deleting the edge~$e$ and do by $G''_e$ the graph obtained by contracting the edge~$e$. Consider propagators~(\ref{Delta}) for some fixed~$e$ which is neither a loop nor an isthmus. Let us divide each of sums~(\ref{defFG}) and~(\ref{defwFG}) into two parts. Let the first part of the propagator of the fixed edge~$e$ contain its constant (which equals $a$ and $\tilde a$, correspondingly). Let the second part contain the remaining term (which equals, correspondingly, $b\,\delta(x_{i(e)}-x_{f(e)})$ and $\tilde b\,\delta(k_e)$). As a result we obtain correlations
\begin{eqnarray}
&{\mathcal F}_G(q;a,b)=a\, {\mathcal F}_{G'_e}(q;a,b)+b\,{\mathcal F}_{G''_e}(q;a,b),&\label{FA1}\\
&\widetilde{\mathcal F}_G(q;\tilde a,\tilde b)=\tilde b\, \widetilde{\mathcal F}_{G'_e}(q;\tilde a,\tilde b)+
\tilde a\,\widetilde{\mathcal F}_{G''_e}(q;\tilde a,\tilde b).& \label{FA2}
\end{eqnarray}
Recurrent correlations in form~(\ref{FA1},\ref{FA2}) are called Tutte--Grothendieck recurrent correlations~\cite{welsh-book}.

Recall that the Whitney rank generating function $R_M(u,v)$ of the matroid $M=(E,r)$ obeys the formula
\begin{equation}
\label{Whitney}
R_M(u,v)=\sum_{A\subseteq E}u^{r(E)-r(A)}v^{|A|-r(A)},
\end{equation}
while the Tutte polynomial is defined as $T_{M}(u,v)=R_{M}(u-1,v-1)$. Evidently that $T_{M^\perp}(u,v)=T_{M}(v,u)$ and
\begin{equation}
\label{chifromR}
\chi_{M}(x)=(-1)^{r(E)}R_{M}(-x,-1),\quad \chi_{M^\perp}(x)=(-1)^{|E|-r(E)}R_{M}(-1,-x).
\end{equation}

In this section we consider only cyclic matroids of graphs $G=(V,E)$. Their rank function for any subgraph $A$ ($A\subseteq E$) of the graph $G$ is defined as $r(A)=|V(A)|-|K(A)|$. Recall~(\cite{tutte-book}) that the dichromatic polynomial of the graph $Q_G(u,v)$ obeys the formula $u^{|K(G)|}R_M(u,v)$, where $M$ is the cyclic matroid of the graph. Therefore $Q_G(u,v)$ obeys the formula
$$
Q_G(u,v)=\sum_{A:A\subseteq E} u^{|K(A)|} v^{|A|-|V|+|K(A)|}.
$$
As appeared, the following exact formulas are valid for FA mentioned in Proposition~\ref{ter:3}:
\begin{eqnarray}
&\mathcal F_G(q;a,b) = a^{|E|-|V|} b^{|V|} Q_G(q\,a/b,\,b/a),& \label{t1}\\
&\widetilde{\mathcal F}_G(q;\tilde a,\tilde b) = \tilde a^{\,|V|}\, \tilde b^{\,|E|-|V|}\,
Q_G (\tilde b/\tilde a,\, q\,\tilde a/\tilde b). &\label{t2}
\end{eqnarray}
and Proposition~\ref{ter:3} evidently follows from these formulas.

Let us finally define notions of a ``bad coloring polynomial'' and a ``bad flow polynomial'' as polynomials of two variables which are usually denoted by~$q$ and~$x$. The notion of a bad coloring polynomial was introduced in the book~\cite[chapter~4]{welsh-book} of Dominic Welsh in connection with the calculation of the partition function in the Potts model;
it is also known as a coboundary polynomial. Its definition coincides with the definition of ${\mathcal F}_G(q;1,x-1)$. The notion of a bad flow polynomial (a boundary polynomial) was introduced analogously in~\cite{goodall}; it coincides with $\widetilde{\mathcal F}_G(q;1,x-1)$.

Bad coloring and bad flow polynomials can be expressed in terms of Tutte polynomials analogous to formulas~(\ref{t1},\ref{t2}); the formula which establishes a connection between them is also known, it is analogous to~(\ref{eq:terfourier}). Since, evidently,
$$
{\mathcal F}_G(q;a,b)=a^{|E(G)|} {\mathcal F}_G(q;1,b/a),\quad
\widetilde{\mathcal F}_G(q;\tilde a,\tilde b)=\tilde a^{|E(G)|} \widetilde{\mathcal F}_G(q;1,\tilde b/\tilde a),
$$
these formulas (stated as theorem~2~in~\cite{goodall}) are equivalent to formulas~(\ref{t1},\ref{t2}) and Proposition~\ref{ter:3}.

\begin{sled} {\rm \textbf{(Theorem~\ref{ter:2} for graphic matroids)}}\
Formulas~(\ref{kuptscov1}),~(\ref{kuptscov2}), and~(\ref{twograph}) are valid.
\end{sled}
\textbf{Proof:}\ By definition, $F_G(q)=\widetilde{\mathcal F}_G(q;1,-1)$. Consequently, by Proposition~\ref{ter:3} we have $F_G(q)=q^{-|V|}\mathcal F_G(q;-1,q)$. Each propagator $\Delta_e(x)=-1+q\delta(x)$ that enter in the right-hand side of this equality is representable as follows:
\begin{equation}
\label{Deltaex1}
\Delta_e(x)=(-q)(1-\delta(x))+(q-1).
\end{equation}
We get
$$
F_G(q)=q^{-|V|} \sum_{A:A\subseteq E} (-q)^{|A|} P_A(q) (q-1)^{|E|-|A|},
$$
which is equivalent to formula~(\ref{kuptscov1}).

In order to immediately deduce~(\ref{kuptscov2}) from~(\ref{eq:terfourier}), note that the equality
$$
P_G(q)=q^{|V|}\widetilde{\mathcal F}_G(q;-1/q,1)
$$
is a particular case of the above correlation. Representing each pro\-pa\-ga\-tor~$\Delta_e(k)=-1/q+\delta(k)$ in the right-hand side of this equality as the sum
$$
\Delta_e(k)=(-1/q)(1-\delta(k))+(1-1/q)\delta(k),
$$
we obtain
$$
P_G(q)=q^{|V|} \sum_{A:A\subseteq E} (-1/q)^{|A|} F_A(q) (1-1/q)^{|E|-|A|},
$$
which is equivalent to formula~(\ref{kuptscov2}).

Finally, in order to immediately deduce~(\ref{twograph}) from~(\ref{eq:terfourier}), it suffices to use in the first paragraph of our proof the representation
$$
\Delta_e(x)=(-1)(1-\delta(x))+(q-1)\delta(x)
$$
for the propagator $\Delta_e(x)=-1+q\delta(x)$ in place of~(\ref{Deltaex1}).$\square$

\section{Proof of Theorem~\ref{ter:2}}

As was mentioned in Introduction, it is not difficult to immediately prove Theorem~\ref{ter:2} by substituting formula~(\ref{chiM}) in both sides of equalities under consideration and by performing some algebraic transformations. However here we intentionally omit the direct proof, because we make an attempt to establish a connection of Theorem~\ref{ter:2} with other known results, within whose frameworks this theorem occurs more naturally.

In the previous section we have considered the case of graphic matroids. In Section~2 we considered their generalization, namely, $\mathbb F_q$-linear matroids. As an intermediate result, when proving Theorem~\ref{ter:1} we obtained Lemma~\ref{lem:1}, as well as Lemma~\ref{utv:1} which is a generalization of the previous one to the case of any natural values of the parameter~$j$ (in Lemma~\ref{lem:1} one considers only cases of $j=2$ and $j=1$). Let us obtain Theorem~\ref{ter:2} as a corollary of any of assertions of the mentioned lemmas with some fixed value of the parameter~$j$.

\begin{sled} {\rm \textbf{(Theorem~\ref{ter:2} for $\mathbb F_q$-linear matroids)}}\
\label{sled:3}
Formulas~(\ref{one}) and~(\ref{two}) are valid for $\mathbb F_q$-linear ma\-tro\-ids.
\end{sled}

We need the following propositions.

\begin{utv}
\label{utv:aboutLinearEquationOverFq}
\quad Fix arbitrary constant values $c_1,\ldots,c_\ell$ in~$\mathbb F_q^*$ and variables ${\alpha_1,\ldots,\alpha_\ell}$ in~$\mathbb F_q^*$. Then any inhomogeneous equation $\sum_{i=1}^k c_i\alpha_i=b$, ${b\in\mathbb F_q^*}$, has one and the same number of solutions; this number is less than the number of solutions of the homogeneous equation $\sum_{i=1}^\ell c_i\alpha_i=0$ by~$(-1)^\ell$.
\end{utv}
\textbf{Proof:} The first part of the proposition is evident. The second one follows from the fact that the desired difference is representable as the sum of values of the additive character~(\ref{simplechi}), i.e., various values of roots of the~$p$th degree of one, i.e.,
$$
\sum_{\alpha\in \left(\mathbb F_q^*\right)^\ell}\chi\left(\sum_{i=1}^\ell c_i\alpha_i\right)=
\prod_{i=1}^\ell \sum_{\alpha_i\in \mathbb F_q^*}\chi\left(c_i\alpha_i\right)=(-1)^\ell .
$$

\begin{utv}
\label{utv:Rota}
Let us use the following denotations: $M=(E,r)$ is an $\mathbb F_q$-linear matroid; ${\,\tt M\,}$ is the matrix of its representation with the set of rows~$V$, $|V|=r(E)$; $A: A\subseteq E$ is a fixed subset of~$E$; ${\mathbf x}=(x_v, v\in V)$ is a row of variables with values in $\mathbb F_q^V$;
${\mathbf y}={\mathbf x}{\,\tt M\,}$. Then $\chi_{M/A}(q)$ coincides with the number of all possible vectors~${\mathbf x}$ such that $y_e=0$, if $e\in A$, and $y_e\neq 0$ otherwise.
\end{utv}
\textbf{Proof:} The case of this proposition, when $A=\emptyset$, is well known; it coincides with a part of the Crapo--Rota theorem related to a famous critical problem. The following considerations reduce the case when the set~$A$ is nonempty to the case of the empty set.

It is well known (see \cite[Proposition 3.2.6]{Oxley}) that the matrix that corres\-ponds to the matroid~$M/e$ can be obtained from the matrix~${\,\tt M\,}$ with the help of one of two operations.
\begin{itemize}
\item If the column~$e$ is zero, then it suffices just to delete it. Note that in this case the element~$y_e$ of the vector~${\mathbf y}$ always equals zero independently of ${\mathbf x}$.
\item If the column~$e$ is nonzero, then we need to perform some elementary transformations of rows of the matrix~${\,\tt M\,}$ which make all elements of the column~$e$, except one of them, namely, $(v,e)$, vanish; it is necessary to delete the row~$v$ and the column~$e$ from the obtained matrix.
Note that elementary transformations of rows of the matrix~${\,\tt M\,}$ are equivalent to bijective transformations in the space conjugate to the space of vectors~${\mathbf x}$,
and the number of vectors~${\mathbf x}$ which turn elements~$y_e$ to zero (or not), $e\in E$, remains constant. By deleting the row~$v$ and the column~$e$ we put $x_v=0$ in terms of new variables. This is the only possible way to make the element~$y_e$ of the vector~${\mathbf y}$, which corresponds to the transformed matrix, turn into zero.
\end{itemize}

These considerations in fact prove Proposition~\ref{utv:Rota} in the case of ${A=\{e\}}$. Since
$$
M/A=(\ldots((M/e_1)/e_2)/\ldots)/e_\ell,\quad \mbox{where $\{e_1,\ldots,e_\ell\}=A$},
$$
by repeating the same considerations~$\ell$ times, we conclude that in order to calculate~$\chi_{M/A}$ we need not to consider the whole vector space~${\mathbf x}$, but only its subspace such that $y_e=0$ for all $e\in A$. For this subpace we can apply the Crapo--Rota theorem. $\square$

\textbf{Proof of Corollary~\ref{sled:3}:} In order to calculate the number of solutions of the equation $Q_j({\mathbf x},\alpha)=b$ we enumerate (for various $A\in E$) numbers of solutions of systems of two equations:\\
1) $y_e=0$ for all $e\in A$; $(y_e)^j=c_e \neq 0$ for all $e\not\in A$,\\
2) $\sum_{e\in E\setminus A} c_e\alpha_e=b$. \\
Note that by Proposition~\ref{utv:aboutLinearEquationOverFq} the number of solutions of the second equation is independent of the choice of concrete values of constants $c_e$. Therefore we require only that $c_e\ne 0$.

According to Proposition~\ref{utv:Rota}, the number of solutions of the first equation with respect to variables~${\mathbf x}$ is $\chi_{M/A}(q)$. By Proposition~\ref{utv:aboutLinearEquationOverFq} the difference of numbers of solutions of the homogeneous and inhomogeneous variants of the second equation in variables $\alpha\in\left(\mathbb F_q^*\right)^E$ is $(q-1)^{|A|}(-1)^{|E|-|A|}$ (here we take into account the fact that no constraint is imposed on variables $\alpha_e$, $e\in A$). By Lemma~\ref{utv:1} we obtain
$$
\chi_{M^\perp}(q)=q^{-|V|} \sum_{A:A\in E} (-1)^{|E|-|A|} (q-1)^{|A|} \chi_{M/A}(q),
$$
which in view of the condition $r(E)=|V|$ is equivalent to~(\ref{two}).

Note that when proving (in Introduction) the equivalence of formulas~(\ref{one}) and~(\ref{two}) for arbitrary matroids we used only two operations, namely, the calculation of minors of matroids and the transition to dual matroids. Since these operations are closed in the class of $\mathbb F_q$-linear matroids, for all matroids in this class we have just proved both formulas~(\ref{two}) and~(\ref{one}).
\noindent$\square$

In earlier papers we have already established equalities, whose left-hand sides contain Tutte polynomials indexed by the matroid~$M$, while the right-hand ones represent the sum over subsets~$A$ of the set~$E$, and the addends are expressed in terms of polynomials indexed by matroids $M|A$ and $M/A$. The first and most known of these formulas, namely, the convolution formula for the Tutte polynomial, has appeared in~\cite{reiner3}. It takes the form
\begin{equation}
\label{rein1}
T_M(x,y)=\sum_{A:A\subseteq E} T_{M|_A} (0,y) T_{M/A} (x,0).
\end{equation}
In the paper ~\cite{reiner2} V.~Reiner generalizes equality~(\ref{rein1}) to a certain two-parameter family of equalities, where formula~(\ref{rein1}) is a limit correlation. In ~\cite{kung1} J.~P.~S.~Kung establishes the following multiplication identity for characteristic polynomials of matroids:
\begin{equation}
\label{kung1}
\chi_M(\lambda\,\xi)=\sum_{A:A\subseteq L(M)} \chi_{M/A}(\lambda)\,\xi^{r(E)-r(M|_A)} \chi_{M|_A}(\xi);
\end{equation}
here the summation is performed over all flats~$L(M)$ of the matroid~$M$. Finally, in~\cite{kung2} one obtains several identities which generalize identities~(\ref{rein1}) and~(\ref{kung1}). In particular, identity~3 in this paper for the Whitney rank ge\-ne\-ra\-ting function~(\ref{Whitney}) takes the form
\begin{equation}
\label{kung2}
R_M(\lambda\xi,xy)=\sum_{A:A\subseteq E} \lambda^{r(E)-r(A)}(-y)^{|A|-r(A)}R_{M|_A}(-\lambda,-x)R_{M/A}(\xi,y).
\end{equation}
Let us represent formulas~(\ref{one}) and~(\ref{two}) as particular cases of formula~(\ref{kung2}).

\begin{utv}
\begin{equation}
\label{last2}
R_M(x,1/x)=(1+1/x)^{|E|}x^{r(E)}.
\end{equation}
\end{utv}
\textbf{Proof:}
$$
R_M(x,1/x)=\sum_{A:A\subseteq E} x^{r(E)-r(A)}x^{r(A)-|A|}=x^{r(E)}\sum_{A:A\subseteq E}x^{-|A|}=
(1+1/x)^{|E|}x^{r(E)}.
$$

\textbf{Proof of Theorem~\ref{ter:2} as a corollary of equality~(\ref{kung2}):} Putting $xy=-q$ and $\lambda\xi=-1$, in accordance with~(\ref{chifromR}) in the left-hand side of equality~(\ref{kung2}) we obtain the necessary expression $(-1)^{|E|-r(E)}\chi_{M^\perp}(q)$. For establishing~(\ref{one}) it suffices to concretize $x=1$, $\lambda=q$ and, correspondingly, $y=-q$, $\xi=-1/q$; then to use substitutions~(\ref{chifromR}) and~(\ref{last2}), and sum degrees of~$q$ and~$(-1)$. Formula~(\ref{two}) can be obtained analogously, provided that $y=-1$, $\xi=-q$ and, correspondingly, $x=q$ and $\lambda=1/q$.

\section{Conclusion}
Formulas for the characteristic polynomial of the dual matroid obtained in this paper were initiated by methods of the theory of FA applied to amplitudes over finite fields. Note that recurrent correlations for vacuum FA~(\ref{FA1},\ref{FA2}) lead to a new ``discovery'' of Tutte polynomials.

The representation of the characteristic polynomial of a {matroid} as a linear combination of Legendre symbols is less known. Such a representation was first obtained for the flow polynomial. Lemma~\ref{lem:1} solves the puzzle of the fact that the linear combination of Legendre symbols with coefficients $\pm \frac1{q^n}$ gives the number of everywhere nonzero flows. As appeared, formula~(\ref{eq:3}) occurs as a result of the re-ordering performed in the calculation of the number of solutions of the equation~$Q({\mathbf x},\alpha)=0$.

Note that the flow polynomial coincides with the vacuum FA. The nonva\-cuum variant of FA over a finite field (see \cite[Section 3]{lerner}) coincides with the number of all possible nonzero flows with fixed values of variables~$k_e$ associated with all edges~$e$ which are incident to some vertex~$v'$ (these edges form the set of outer edges~$E_{ext}$, and the corresponding variables are called \textit{external momenta}, the nonvacuum Feynman amplitude depends on their values). Let us impose the following natural constraint on external momenta:
\begin{equation}
\label{eq:conc}
\sum_{e\in E_{ext}} \varepsilon_{v'e} k_e=0.
\end{equation}
It is interesting when a nonvacuum FA is positive with all nonzero external momenta satisfying condition~(\ref{eq:conc}) (in what follows, for brevity, when speaking about the positiveness of a nonvacuum FA, we assume that this condition is fulfilled).

In the paper~\cite{kochol} Martin Kochol proves (using the notion of ``strongly $\mathbb Z_q$-connected'' FA instead of nonzero nonvacuum ones) that the Tutte 3-flow conjecture is equivalent to the positiveness of nonvacuum FA over~$\mathbb F_3$ for any 4-edge-connected multigraph~$G$ with vertices of the degree~4 or~5 with any choice of the vertex~$v'$ (and values of external momenta~$k_e$, $e\in E_{ext}$). One can easily prove that the Tutte 5-flow conjecture is equivalent to the positiveness of nonvacuum amplitudes over~$\mathbb F_5$ for any cubic three-connected graphs~$G$; as above, the vertex~$v'$ and values of external momenta can be chosen arbitrarily.

Initially the theory of FA was developed for a nonvacuum case. Evidently, there exist nonvacuum analogs of theorems~\ref{ter:1} and~\ref{ter:2} which allow one to calculate the number of corresponding flows. However due to the complexity of corresponding nonvacuum formulas they can hardly be efficiently used in practice. Nevertheless, in view of remarks made in the previous paragraph, both the notion of nonvacuum FA and the traditional application of these objects in the theory of FA, possibly, will be useful in combinatorics.

\section*{Acknowledgements} The author is grateful to A.P.Kuptsov, in collaboration with whom in 2002--2003 he obtained all basic results described in Section~3.


\begin{thebibliography}{10}
\bibitem{aigner}
M.~Aigner, {\it Combinatorial theory}, Grundlehren der Mathematischen Wissenschaften, 234, Springer, Berlin, 1979. \MR{0542445}

\bibitem{belk}
P.~Belkale\ and\ P.~Brosnan, Matroids, motives, and a conjecture of Kontsevich, Duke Math. J. {\bf 116} (2003), no.~1, 147--188. \MR{1950482}

\bibitem{chung}
F.~Chung\ and\ C.~Yang, On polynomials of spanning trees, Ann. Comb. {\bf 4} (2000), no.~1, 13--25. \MR{1763947}

\bibitem{goodall}
A.~J.~Goodall, Some new evaluations of the Tutte polynomial, J. Combin. Theory Ser. B {\bf 96} (2006), no.~2, 207--224. \MR{2208351}

\bibitem{kochol}
M.~Kochol, An equivalent version of the 3-flow conjecture, J. Combin. Theory Ser. B {\bf 83} (2001), no.~2, 258--261. \MR{1866721}

\bibitem{reiner3}
W.~Kook, V.~Reiner\ and\ D.~Stanton, A convolution formula for the Tutte polynomial, J. Combin. Theory Ser. B {\bf 76} (1999), no.~2, 297--300. \MR{1699230}

\bibitem{kung1}
J.~P.~S.~Kung, A multiplication identity for characteristic polynomials of matroids, Adv. in Appl. Math. {\bf 32} (2004), no.~1-2, 319--326. \MR{2037633}

\bibitem{kung2}
J.~P.~S.~Kung, Convolution-multiplication identities for Tutte poly\-no\-mials of graphs and matroids, J. Combin. Theory Ser. B {\bf 100} (2010), no.~6, 617--624. \MR{2718681}



\bibitem{lerner}
E.~Yu.~Lerner, Feynman integrals of $p$-adic argument in the momentum space. II. Explicit expressions,
Teoret. Mat. Fiz. {\bf 104} (1995), no. 3, 371--392; translation in Theoret. and Math. Phys. {\bf 104} (1995), no.~3, 1061--1077 (1996). \MR{1607005}

\bibitem{ourArX}
E.~Yu.~Lerner, A.~P.~Kuptsov, S.~A.~Mukhamedjanova, Flow polynomials as Feynman amplitudes and their $\alpha$-representation,
\href{https://arxiv.org/abs/1609.01120}{arXiv:1609.01120}, 2016, submitted to Electronic Journal of Combinatorics.

\bibitem{lidl}
R.~Lidl\ and\ H~Niederreiter, {\it Finite fields}, second edition, Encyclopedia of Mathematics and its Applications, 20, Cambridge Univ. Press, Cambridge, 1997. \MR{1429394}

\bibitem{Oxley}
J.~G.~Oxley, {\it Matroid Theory}, Oxford University Press, Oxford, New York, 1992. \MR{1207587}

\bibitem{reiner2}
V. Reiner, {\it An interpretation for the Tutte polynomial}, European J. Combin. {\bf 20}:2
(1999), 149--161. \MR{1676189}

\bibitem{reiner1}
V. Reiner, Lectures on Matroids and Oriented Matroids, 2005. Available
\href{https://www.math.umn.edu/reiner/Talks/Vienna05/Lectures.pdf}
{www.math.umn.edu/reiner/Talks/Vienna05/Lectures.pdf}.

\bibitem{schmidt}
W.~M.~Schmidt, {\it Equations over finite fields. An elementary approach}, Lecture Notes in Mathematics, Vol. 536, Springer, Berlin, 1976. \MR{0429733}

\bibitem{stanleyArticle}
R.~P.~Stanley, Spanning trees and a conjecture of Kontsevich, Ann. Comb. {\bf 2} (1998), no.~4, 351--363. \MR{1774974}

\bibitem{tutte-book}
W.~T.~Tutte, {\it Graph theory}, reprint of the 1984 original, Encyclopedia of Mathematics and its Applications, 21, Cambridge Univ. Press, Cambridge, 2001. \MR{1813436}

\bibitem{welsh-book}
D.~J.~A.~Welsh, {\it Complexity: knots, colourings and counting}, London Mathematical Society Lecture Note Series, 186, Cambridge Univ. Press, Cambridge, 1993. \MR{1245272}


\end{thebibliography}
\end{document}